\documentclass[12pt]{article} \usepackage{amssymb,amsmath,amscd}

\newtheorem{proposition}{Proposition}
\newtheorem{corollary}[proposition]{Corollary}
\newtheorem{lemma}[proposition]{Lemma}
\newtheorem{remark}[proposition]{Remark}

\def\Co{\mathbb{C}} 

\def\fidi{\hskip5pt \vrule height4pt width4pt depth0pt \par}

\def\ker{{\rm Ker}}

\def\Na{\mathbb{ N}}  
\def\Pol{{\rm Pol}}

\def\Re{\mathbb{R}}  

\def\sphere{\mathbb{S}}
\def\sph4{\Sigma_q^4}
\def\sph2n{\Sigma_q^{2n}}
\def\Tr{{\rm Tr}}
\def\tr{{\rm tr}}
\def\u4{{\bf u}(4)}

\def\Ze{{\mathbb Z}}

\def\M{{\cal M}}

\begin{document}

\title{Quantum even spheres $\Sigma^{2n}_q$\\ from
Poisson double suspension}

\author{F.Bonechi${}^{1,2}$, N.Ciccoli${}^3$, M.Tarlini${}^{1,2}$}

\date{\today}

\maketitle

\centerline{{\small  ${ }^1$ INFN Sezione di Firenze}}
\centerline{{\small ${ }^2$ Dipartimento di Fisica, Universit\`a
di Firenze, Italy. }} \centerline{{\small ${ }^3$ Dipartimento di
Matematica, Universit\`a di Perugia, Italy. }} \centerline{{\small
e-mail: bonechi@fi.infn.it, ciccoli@dipmat.unipg.it,
tarlini@fi.infn.it}}
\begin{abstract}
\noindent
{We define even dimensional quantum spheres $\Sigma^{2n}_q$ that
generalize to higher dimension the {\it standard} quantum
two-sphere of Podle\'s and the four-sphere $\Sigma^4_q$ obtained
in the quantization of the Hopf bundle. The construction relies on
an iterated Poisson double suspension of the standard Podle\'s
two-sphere. The Poisson spheres that we get have the same
symplectic foliation consisting of a degenerate point and a
symplectic plane and, after quantization, have the same
$C^*$--algebraic completion. We investigate their $K$-homology and
$K$-theory by introducing Fredholm modules and projectors.}
\end{abstract}

\medskip
\noindent
{\bf Math.Subj.Classification}: 14D21, 19D55, 58B32, 58B34, 81R50


\thispagestyle{empty}

\bigskip
\bigskip

\section{Introduction}
In the seminal paper \cite{Podles} by Podle\'s on
$SU_q(2)$-covariant quantum two-spheres it was introduced a family
of deformations $\sphere^2_{q,d}\ $, depending on $d\in\Re$ and
$q<1$. They define three distinct quantum topological spaces
$\sphere^2_{q,d}\ $: the case $d=0$, the so called {\it standard
sphere}, $d>0$, the {\it non standard sphere}, and $d=-q^{2n},\
n=1,2\ldots$, the {\it exceptional case}.

The geometry of these quantum spaces can be nicely interpreted by
looking at the underlying Poisson geometry and considering the
sphere as a patching of symplectic leaves. There exists, in fact,
a $1$--parameter family of $SU(2)$--covariant Poisson bivectors on
the sphere $\sphere^2$ such that $\sphere^2_{q,d}$ can be seen as
a quantization of such structures (\cite{Sh}). In the case $d=0$
the symplectic foliation is made of a point and a symplectic
$\Re^2$; after quantization, the symplectic plane is quantized to
$K$, the algebra of compact operators, and the degenerate point
survives as a character. The $C^*$-algebra $C(\sphere^2_{q,0})$ is
then isomorphic to the minimal unitization of compacts
$\tilde{K}$, and satisfies the following exact sequence
$0\rightarrow K\rightarrow
C(\sphere^2_{q,0})\rightarrow\Co\rightarrow 0$. In the case $d>0$
symplectic foliation consists of an $\sphere^1$-family of
degenerate points, and two symplectic disks. After quantization
$C(\sphere^2_{q,d})$ satisfies $0\rightarrow K\oplus K\rightarrow
C(\sphere^2_{q,d})\rightarrow C(\sphere^1)\rightarrow 0$. The
exceptional spheres correspond to the symplectic case and after
quantization $C(\sphere^2_{q,-q^{2n}})$ is isomorphic to the
finite dimensional algebra $M_n(\Co)$ of matrices.

In higher dimension, there exist the so called {\it euclidean}
spheres $\sphere^{n}_q$ introduced in \cite{FRT} as quantum
homogeneous spaces of $SO_q(n+1)$; let us notice that in odd
dimension they coincide with the so called Vaksman-Soibelman
spheres $\sphere^{2n+1}_q$ (\cite{VS}) introduced as quantum
quotients $U_q(n)/U_q(n-1)$, as observed in \cite{HL}.

This family of spheres represents a generalization of the Podle\'s
$d=1$ case. In fact they satisfy the following exact sequences
(see \cite{HL,HS,VS})
$$
0\rightarrow K\oplus K\rightarrow C(\sphere^{2n}_q)\rightarrow
C(\sphere^{2n-1}_q)\rightarrow 0
$$
$$
0\rightarrow C(\sphere^1)\otimes K\rightarrow C(\sphere^{2n+1}_q)\rightarrow
C(\sphere^{2n-1}_q)\rightarrow 0\, .
$$
This behaviour reflects exactly the Poisson level where each
$\sphere^{2n}$ contains an equator $\sphere^{2n-1}$ as a Poisson
submanifold and the two remaining hemispheres are open symplectic
leaves of dimension $2n$. In particular all these spheres contain
an $\sphere^1$-family of degenerate points.

As there are no symplectic forms on higher dimensional spheres,
exceptional spheres are possible only in dimension $2$.

In this paper we introduce the family of even spheres
$\Sigma^{2n}_q$ which generalize the $d=0$ Podle\'s sphere: they
all share with it the same kind of symplectic foliation. The idea
of the construction relies on two classical subjects: double
suspension and Poisson geometry.

The suspension idea is certainly not new and, in fact, appears in
many of the papers devoted to the construction of particular
deformations of the four-sphere (\cite{CoLa,BG,DLM,Si}).  When
one considers double suspension, however, more interesting
possibilities appear: on one hand one could define a \emph{purely
classical} double suspension already at an algebraic level by
adding a pair of central selfadjoint generators and modding out
suitable relations. A different kind of double suspension was
considered, at the $C^*$--algebra level, in \cite{HS}. There the
authors consider the non reduced double suspension of a
$C^*$--algebra $A$ as the middle term $S^2A$ of a short exact
sequence
$$
0\to A\otimes K\to S^2A \to {\cal C}({\sphere^1})\to 0
$$
for a suitable fixed Busby invariant. Let us remark that such a
non reduced double suspension always has a naturally defined
${\sphere^1}$--family of characters. All euclidean spheres were
reconstructed in this way starting either from a two--point space or
from $\sphere^1$.

In this paper we will consider the \emph{reduced} double
suspension or reduced topological product $S^2X=\sphere^2\times
X/\sphere^2\vee X$ where, given $p\in \sphere^2$ and $q\in X$,
$\sphere^2\vee X=(\sphere^2\times q)\cup(p\times X)$.

It is quite natural to look at the interaction between the double
suspension and Poisson geometry. One word has to be said about the
fact that while suspension is an essentially topological
construction, Poisson bivectors are of differential nature, so
that, in principle, on the double suspension of a given manifold
there's no reason to have a manifold structure, let aside a
Poisson bracket. Still whenever the manifold structure is there
one can ask whether such a Poisson structure arises. More
precisely the double suspension of a manifold $M$ can be seen as a
topological quotient $S^2M$ of $\sphere^2\times M$. If we are
given a Poisson bivector on the $2$--sphere and a Poisson bivector
on $M$ we can ask whether the quotient map $\sphere^2\times M\to
S^2M$ coinduces a Poisson bracket on the quotient. If this is the
case we then look for its quantization.

From this point of view the \emph{classical} double suspension
quantizes a double suspension with respect to the trivial Poisson
structure on $\sphere^2$ while Hong--Szyma\'nski construction
corresponds to the standard symplectic structure on $\Re^2$
attached along an $\sphere^1$ which will survive as a family of
$0$--leaves on the suspension.

In this paper the double suspension is built by assuming the
Podle\'s $d=0$ Poisson structure on the two-sphere and considering
$\sphere^{2n}$ as the double suspension of $\sphere^{2(n-1)}$. We
will show that a \emph{Poisson double suspension} of spheres
exists for each $n$ and that it quantizes both at an algebraic
level and at the $C^*$--algebra level. The spheres $\Sigma^{2n}_q$
that we obtain have all the same symplectic foliation of the
two-sphere, and their quantizations are topologically equivalent,
{\it i.e.} they are the minimal unitization of compacts and
satisfy
$$ 0\rightarrow K\rightarrow C(\Sigma^{2n}_q)\rightarrow
\Co\rightarrow 0 ~~. $$

We then have a quite extreme case of {\it quantum degeneracy}:
these quantum spaces, whatever is the classical dimension, are all
topologically equivalent to a zero dimensional compact quantum
space. This is an extreme manifestation of a well known fact, that
quantum spaces associated to quantum groups have lower dimension
than the classical one. Moreover this reminds of canonical
quantization in which the Weyl quantization of $C_0(\Re^{2n})$ is $K$,
for each $n$. The opposite behavior is
represented by the so called $\theta$-deformation, whose behavior
is almost classical (see \cite{CoLa,CoDV}).

In the case of the four-sphere $\Sigma^4_q$, the algebra that we
get is that obtained in \cite{BCT1,BCT2}, in the context of a
quantum group analogue of the Hopf principal bundle
$\sphere^7\to\sphere^4$. It is unclear whether all even spheres
$\Sigma^{2n}_q$ can be obtained as coinvariant subalgebras in
quantum groups.

In Section 2 we introduce the Poisson double suspension that
iteratively defines the Poisson even spheres and we study their
symplectic foliation. In Section 3 we introduce the quantization
$\Pol(\Sigma^{2n}_q)$ at the level of polynomial functions, we
classify the irreducible representations in bounded operators and
then show that the universal $C^*$-algebra $C(\Sigma_q^{2n})$ is
the minimal unitization of compacts. We introduce Fredholm modules
for each of these spheres. In Section 4, we give the non trivial
generator of $K_0(\Sigma^{2n}_q)$ and compute its coupling with
the character of the previously introduced Fredholm modules.

\bigskip
\bigskip
\section{The standard Poisson structure.}

Let us define a point on $\sphere^2\times\ldots\times \sphere^2$
by giving to it coordinates $((\alpha_1,\tau_1),\ldots
(\alpha_n,\tau_n))$, where $|\alpha_i|^2=\tau_i(1-\tau_i)$. Let
$M$ be the matrix whose entries are $M_{ii}=0$, $M_{ij}=1$ and
$M_{ji}=1/2$ if $i<j$. Let \begin{equation} \label{mappa} a_i=
\alpha_i \prod_k \tau_k^{M_{ik}}\quad t=\prod_i \tau_i\;.
\end{equation}

Since \begin{eqnarray*} \sum_i |a_i|^2 &=& \sum_i
\tau_i(1-\tau_i)\prod_k \tau_k^{2M_{ik}} =
\tau_1(1-\tau_1)\tau_2^2\ldots \tau_n^2 +\cr & &
\tau_2(1-\tau_2)\tau_1\tau_3^2\ldots \tau_n^2+ \ldots +
\tau_n(1-\tau_n) \tau_1\ldots \tau_{n-1}\cr &=& t(1-t)\;,
\end{eqnarray*} then relation (\ref{mappa}) defines a projection
into $\sphere^{2n}$. Let us call $\Phi:
\sphere^2\times\ldots\times \sphere^2\rightarrow \sphere^{2n}$
such projection. One can verify that this map is equivalent to the
iterated reduced double suspension of a two-sphere, with preferred
point the North Pole $\alpha=\tau=0$. In fact the map from the
cartesian product $\sphere^2\times X$ to the reduced topological
product $S^2X$ is the unique continuous map which is
a homeomorphism everywhere but on the counter image of a point
over which its fiber is $\sphere^2\vee X$. Starting from the
two-sphere and iterating this procedure one defines a map from
$\sphere^2\times\ldots\sphere^2$ to $\sphere^{2n}$ that is a
homeomorphism everywhere but on the North Pole, where its fiber is
the topological join of $n$ copies of
$\underbrace{\sphere^2\times\ldots\times\sphere^2}_{n-1}$. This map is the
projection $\Phi$.

\medskip

Let us equip $\sphere^2$ with the standard Poisson structure, {\it
i.e.} the limit structure of the Podle\'s standard two sphere
(\cite{Podles,Sh}), and $\sphere^2\times\ldots \times\sphere^2$
with the product Poisson structure. The brackets among polynomial
functions are defined by giving \begin{equation}\label{poissonS2}
\{\alpha_i,\tau_j\} = -2\delta_{ij}\alpha_i \tau_i\,,~~~~~~
\{\alpha_i,\alpha^*_j\}=2\delta_{ij}(\tau_i^2-\alpha_i\alpha_i^*)
\;. \end{equation} We prove the following result.

\medskip \begin{proposition} \label{Phi} The map $\Phi$ is a Poisson
map. The coinduced brackets on $\sphere^{2n}$ read:
\begin{eqnarray}\label{poissonS4}
\{a_k,a_\ell\} = a_k a_\ell \ \ (k<\ell)\,,&&
\{a_k,a_\ell^*\}=-3a_ka_\ell^*\ \ (k\not=\ell)\,,\cr
\{a_i,\,t\,\}\, =\,-2 a_i t\,,\hspace{1.25cm} && \{a_k,a_k^*\}= 2t^2 +2
\sum_{\ell < k} a_\ell a_\ell^*- 2 a_k a_k^*\;.
\end{eqnarray}
\end{proposition}
{\it Proof}. The result is obtained by explicitly computing the
brackets on $\sphere^2\times\ldots\times\sphere^2$. The only relation
that deserves some attention is the last one. By direct
computation we obtain $$ \{a_k,a_k^*\}= 2
\tau_k^2\prod_i\tau_i^{2M_{ki}} -2 a_ka_k^* \;. $$ Let us show by
induction on $k$ that $$ \tau_k^2\prod_i\tau_i^{2M_{ki}} = t^2 +
\sum_{\ell<k} a_\ell a_\ell^*\;. $$ It is clearly true for $k=1$.
Let it be true for $k$. We then have
\begin{eqnarray*}
\tau_{k+1}^2\prod_i\tau_i^{2M_{{k+1}i}}&=&\tau_1\ldots \tau_k
\tau_{k+1}^2 \ldots \tau_n^2=\tau_1\ldots
(\tau_k^2+\alpha_k\alpha_k^* )\tau_{k+1}^2\ldots \tau_n^2\cr &=&
t^2 + \sum_{\ell<k} a_\ell a_\ell^* + a_k a_k^* = t^2 +
\sum_{\ell<k+1} a_\ell a_\ell^*\;. \end{eqnarray*} \fidi
\medskip\medskip

\noindent The Poisson manifold ($\sphere^{2n},\{,\}$) is
usually called coinduced from the bracket on
$\sphere^2\times\ldots\times\sphere^2$ (see \cite{Vai}). Its
symplectic foliation is described in the following proposition.

\medskip
\begin{proposition}
There are two distinct symplectic leaves in
$(\sphere^{2n},\{,\})$:
\begin{itemize}
\item[$i)$] a zero dimensional leaf given by the north pole
$P_N$=$(a_i=0,t=0)$; \item[$ii)$]
$\Re^{2n}=\sphere^{2n}\setminus{P_N}$.
\end{itemize}

The Poisson brackets on $\Re^{2n}$ read:
\begin{eqnarray}
\label{poisson_plane} \{z_k,z_\ell\} = z_k z_\ell \ \ (k\leq \ell)\,,&&
\{z_k,z_k^*\}= 2(1+\sum_{\ell\leq k}z_\ell
z_\ell^*)\,,\cr
\{z_k,z_\ell^*\}= z_k z_\ell^* \ \ (k\not=\ell)\,.&&
\end{eqnarray}
\end{proposition}

{\it Proof.} It is clear that the north pole $P_N$ defined by
$a_i=t=0$ is a degenerate point. We are going to show that
$\Re^{2n}=\sphere^{2n} \setminus P_N$ is symplectic. Relations
(\ref{poisson_plane}) are obtained by direct computation of the
brackets among the complex coordinates $z_i=a_i/t$.

Let us define the $2n\times 2n$ antisymmetric matrix $S^{(n)}$ as
$S^{(n)}_{ij}=\{w_i,w_j\}$, where $w_{2k-1} = z_k$ and
$w_{2k}=z_k^*$, for $k=1,\ldots,n$. It is clear that
$S^{(n)}_{ij}=S^{(n-1)}_{ij}$ for $i,j=1\ldots,2(n-1)$. To compute
the determinant of this matrix let us introduce a set of $2n$
fermionic variables $\eta_i$, {\it i.e.}
$\eta_i\eta_j+\eta_j\eta_i=0$. The pfaffian of $S^{(n)}$ can be
expressed as
\begin{eqnarray*}
Pf(S^{(n)}) &=& \int d\eta_{2n}\ldots d\eta_1 \
e^{\frac{1}{2}\sum_{ij}S^{(n)}_{ij}\eta_i\eta_j} \cr
&=& \int d\eta_{2n}\ldots d\eta_1 \
e^{\frac{1}{2}\sum_{ij}^{2(n-1)}S^{(n)}_{ij}\eta_i\eta_j}
e^{\sum_{i=1}^{2(n-1)} \eta_i J_i}
e^{S^{(n)}_{2n-1,2n}\eta_{2n-1}\eta_{2n}} \;, \cr
\end{eqnarray*}
where $J_i= S^{(n)}_{i,2n-1}\eta_{2n-1}+S^{(n)}_{i,2n}\eta_{2n}$.
By using standard rules for fermionic integration (see for
instance \cite{Ra}) and taking into account that $J_i J_j=0$ we
get
\begin{eqnarray*}
Pf(S^{(n)}) &=& Pf(S^{(n-1)})\int d\eta_{2n} d\eta_{2n-1}
(1+S^{(n)}_{2n-1,2n}\eta_{2n-1}\eta_{2n})\cr
&=& Pf(S^{(n-1)})S^{(n)}_{2n-1,2n} =
Pf(S^{(n-1)})\{z_n,z_n^*\}\cr
&=& 2 Pf(S^{(n-1)}) (1+\sum_{\ell\leq n}|z_\ell|^2)\;.
\end{eqnarray*}
Since $Pf(S^{(1)})=2(1+|z|^2)$ we conclude that
$Pf(S^{(n)})\not=0$ and that $\Re^{2n}$ is symplectic. \fidi

\medskip

\noindent The Poisson structure on the other chart
$\Re^{2n}=\sphere^{2n}\setminus\{t=1\}$ is symplectic everywhere
but the origin. In \cite{Zak} Zakrzewski introduced a family of
$SU(n)$-covariant Poisson structures on $\Re^{2n}$ with this
foliation. It is an interesting problem to understand the
relations between them.

\bigskip
\bigskip
\section{Quantization of the standard Poisson structure.}

The algebra $\Pol(\sphere^2_{q,0})$  of the Podle\'s standard
sphere is generated by $\{\alpha,\alpha^*,\tau\}$ ,where $\tau$ is
real, with the following relations:
$$
\alpha \tau = q^2 \tau \alpha\ ,\quad q^2 \alpha^* \alpha= \tau (1-\tau),
\quad \alpha\alpha^*=q^2 \alpha^*\alpha+(1-q^2)\tau^2\;.
$$
Let $0 < q < 1$. There are two irreducible representations of
$\Pol(\sphere^2_{q,0})$ with bounded operators, the first is one
dimensional $\epsilon(\alpha)= \epsilon(\tau)=0$, the second
$\sigma:\Pol(\sphere^2_{q,0})\to B(\ell^2(\Na))$ is defined by
\begin{eqnarray}
\sigma(\alpha)|n\rangle &=&
q^{n-1}(1-q^{2n})^{1/2}|n-1\rangle\;,\cr
\label{irrep} \sigma(\tau)|n\rangle &=& q^{2n}
|n\rangle\;.
\end{eqnarray}
Let us define $\sigma^{\otimes n}:\Pol(\sphere^2_{q,0})^{\otimes
n}\to B(\ell^2(\Na)^{\otimes n})$ as the $n^{th}$--tensor product
of the representation $\sigma$, we denote by
$\{\alpha_i,\alpha_i^*,\tau_i\}$ the generators of the
$i^{\,th}\;$--$\;\Pol(\sphere^2_{q,0})$.

\medskip \begin{proposition} \label{sph2n} Let us define
$\Pol(\sph2n)$ the algebra generated by $\{a_i,a_i^*,t\}$ with
relations
$$ a_i\, t=q^2 t\, a_i\;,\quad a_i\, a_j=q^{-1} a_j\,
a_i\;, \quad a_i\, a_j^*=q^{3} a_j^*\, a_i\quad (i < j)\,, $$ $$
a_i\, a_i^*=q^2a_i^*\, a_i+q^2(1-q^2)\sum_{\ell<i}a_\ell^*\,
a_\ell + (1-q^2)t^2\;,\quad \sum_{i=1}^n q^2 a_i^*\, a_i=t-t^2\;.
$$ The mapping $\sigma_n:\Pol(\sph2n) \to B(\ell^2(\Na)^{\otimes
n})$ given by \begin{equation} \label{mappa_q}
\sigma_n(a_i)=\sigma^{\otimes n}(\alpha_i \prod_k
\tau_k^{M_{ik}})\;,\quad \sigma_n(t)=\sigma^{\otimes n}(\prod_i
\tau_i)\;, \end{equation} is a representation of $\Pol(\sph2n)$.
\end{proposition}
\medskip
\noindent {\it Proof.} From the
relation $\sigma^{\otimes n}(\alpha_k \tau_k^{M_{ik}})=
q^{M_{ik}}\sigma^{\otimes n}(\tau_k^{M_{ik}} \alpha_k)$ where
$M_{ij}$ is the matrix with $M_{ii}=0$, $M_{ij}=1$ and
$M_{ji}=1/2$ for $i<j$, it is straightforward to verify the first
line of relations. In order to prove the relations in the second
line we need the following equality (we will omit the application
of $\sigma^{\otimes n}$): $$ \tau_i^2\prod_k
\tau_k^{2M_{ik}}=t^2+q^2\sum_{\ell<i} a_\ell^*\,a_\ell\;. $$ For
$i=1$ it is true, we will verify that it is true for $i+1$
assuming it true for $i$: \begin{eqnarray*} \tau_{i+1}^2\prod_k
\tau_k^{2M_{i+1 k}}&=&\tau_1\cdots\tau_i\tau_{i+1}^2\cdots
\tau_n^2\cr
&=&\tau_1\cdots\tau_{i-1}(q^2\alpha_i^*\alpha_i+\tau_i^2)\tau_{i+1}^2
\cdots \tau_n^2\cr &=&q^2 a_i^*\,a_i +\tau_i^2\prod_k
\tau_k^{2M_{ik}} =q^2 a_i^*\,a_i + t^2+q^2\sum_{\ell<i}
a_\ell^*\,a_\ell\cr &=&t^2+q^2\sum_{\ell<i+1} a_\ell^*\,a_\ell\;.
\end{eqnarray*}
To verify the modulus relation we simply needs the same computation
of the classical case presented at the beginning of Section 2.
\fidi

\bigskip
\noindent
\begin{remark} {\rm The semiclassical limit, defined by
$\{f,g\}=\lim_{q\to 1} \frac{1}{1-q}\,[f,g]$, of the relations of
Proposition \ref{sph2n} coincides with the Poisson structure
defined by the map $\Phi$ in Proposition
\ref{Phi}.\fidi}\end{remark}

\medskip

\medskip
\begin{proposition}
\label{rep_sph2n}
If $\varphi:\Pol(\sph2n) \to  B({\cal H})$ is an irreducible representation
on some Hilbert space, then $\varphi=\epsilon$ or $\varphi=\sigma_n$.
\end{proposition}
\medskip
\noindent
{\it Proof.} In order to prove the existence of an eigenvector of
$\varphi(t)$ we will adapt the proof of Theorem 4.5 in \cite{HMS}.
By using relations we see that $\varphi(t-t^2)>0$ so that
$Sp(\varphi(t))\subset [0,1]$. If $Sp(\varphi(t))=\{0\}$ then
$\varphi(t)=0$ and $\varphi=\epsilon$; if $Sp(\varphi(t)) = \{1\}$
then $\varphi(t)=1$ and this contradicts relations; if
$Sp(\varphi(t)) = \{0,1\}$ then $\lambda = 0$ would be an
eigenvalue and $\ker(\varphi(t))$ would be an invariant subspace.
So in order to have $\varphi$ irreducible and
$\varphi\not=\epsilon$ we must have
$Sp(\varphi(t))\setminus\{0,1\}\not=\emptyset$.

Let $\lambda\in Sp(\varphi(t))\setminus\{0,1\}$ and let
$\{\xi_s\}$ be a set of approximate unit eigenvectors, {\it i.e.}
unit vectors such that
$\lim_{s\rightarrow\infty}||\varphi(t)\xi_s-\lambda\xi_s||=0$. By
writing $t-t^2=\lambda(1-\lambda)+(t-\lambda)(1-\lambda-t)$ we get
\begin{eqnarray*}
||\varphi(t-t^2)\xi_s ||&\geq&
|\lambda(1-\lambda)|-||\varphi(t-\lambda)\varphi(1-\lambda-t)\xi_s||\cr
&\geq&|\lambda(1-\lambda)|-||\varphi(t-\lambda)\xi_s||\
||\varphi(1-\lambda-t)||\geq C' |\lambda(1-\lambda)|\;,
\end{eqnarray*}
for $s$ bigger than some $s_o$ and for some $C'>0$. Moreover we
have
\begin{eqnarray*}
||\sum_{k=1}^n \varphi(a_k^*a_k)\xi_s|| \leq
\sum_{k=1}^n||\varphi(a_k^*)||\ ||\varphi(a_k)\xi_s|| \leq C'' n
||\varphi(a_{k(s)})\xi_s||\;,
\end{eqnarray*}
where $C''$ and $k(s)$ are such that $||\varphi(a_k^*)||\leq C''$
and $||\varphi(a_k)\xi_s||\leq ||\varphi(a_{k(s)})\xi_s ||$ for
all $k$. We conclude that for each $s>s_o$ there exists $1\leq
k(s)\leq n$ such that $||\varphi(a_{k(s)})\xi_s ||\geq
C'''|\lambda(1-\lambda)|$. We can define
$\nu_s=\varphi(a_{k(s)})\xi_s/||\varphi(a_{k(s)})\xi_s||$ and
verify that they are approximating unit eigenvectors for
$q^{-2}\lambda$; in fact
\begin{eqnarray*}
||(\varphi(t)-q^{-2}\lambda)\nu_s || &=&
\frac{q^{-2}}{||\varphi(a_{k(s)})\xi_s
||}||\varphi(a_{k(s)}(t-\lambda))\xi_s ||\cr
&\leq&\frac{q^{-2}\lambda(1-\lambda)}{C'''}||\varphi(a_{k(s)})||\
||\varphi(t-\lambda)\xi_s||\;.
\end{eqnarray*}
We then showed that if $\lambda\in Sp(\varphi(t))\setminus\{0,1\}$
then $q^{-2}\lambda\in Sp(\varphi(t))$. In order to keep
$Sp(\varphi(t))$ bounded it is necessary that for each $\lambda$
there exists $k$ such that $q^{-2k}\lambda = 1$, {\it i.e.}
$Sp(\varphi(t))\setminus\{0\} = \{q^{2k}, k\in \Na\}$. Since each
$q^{2k}$ is isolated, we conclude that it is an eigenvalue.

Let $\psi$ the eigenvector corresponding to $\lambda=1$: since it
is the biggest eigenvalue we get that $\varphi(a_i)\psi=0$. By
direct computation one recovers:
\begin{eqnarray}\label{formula}
\varphi(a_i)\stackrel{\rightarrow}{\prod}_j \varphi(a_j^{*})^{
m_j}\,\psi&=& q^{3\sum_{j<i}m_j}q^{4\sum_{j>i}m_j} q^{2(m_i-1)}
(1-q^{2 m_i})\cr
&&\varphi(\stackrel{\rightarrow} {\prod}_{j<i}a_j^{*\; m_j}\;
a_i^{*\; (m_i-1)}\stackrel{\rightarrow} {\prod}_{j>i} a_j^{*\;
m_j})\,\psi\;.
\end{eqnarray}
Let us define, with ${\bf m}=(m_1,\ldots,m_n)$,
$$\psi^{\bf m}=C^{\bf m}\,
\varphi(\stackrel{\rightarrow} {\prod}_i a_i^{*\; m_i})\psi\;,~~
C^{\bf m}=q^{-(\sum_i m_i)^2+\sum_i m_i(m_i+1)/2}\
\prod_i(q^2;q^2)_{m_i}^{-1/2}\;,$$ where
$(\alpha;q)_s=\prod^s_{j=1}(1-q^{j-1}\, \alpha)$ and
$(\alpha;q)_0=1$. Formula (\ref{formula}) implies
$$||\varphi(a_i^*)\psi^{\bf m}||^2= q^{2\sum_{j\leq
i}m_j}q^{4\sum_{j>i}m_j}(1-q^{2(m_i+1)})||\psi^{\bf m}||^2~~,$$
from which we conclude that $\psi^{\bf m}\not =0$ for each $\bf
m$. The space generated by $\{\psi^{\bf m}\}$ is invariant and it
coincides with ${\cal H}$. Finally it can be verified that the
mapping $T:\ell^2(\Na)^{\otimes n}\rightarrow {\cal H}$ defined by
$|m_1,\cdots,m_n\rangle\,\mapsto\, \psi^{\bf m}$ intertwines the
representations $\sigma_n$ and $\varphi$. Since one can verify
that the $\psi^{\bf m}$'s are orthonormal, we conclude that $T$ is
unitary. \fidi

\medskip
The universal $C^*$-algebra generated by $\Pol(\Sigma^{2n}_q)$ is
then the norm closure of $\sigma_n(\Pol(\Sigma^{2n}_q))$. Since
$\sigma_n(a_i)$ and $\sigma_n(t)$ are trace-class operators, then
$\sigma_n(\Pol(\Sigma^{2n}_q)\setminus\Co)\subset K$. By using
Proposition 15.16 of \cite{DF}, which states that a norm-closed
$*$-subalgebra $A$ of $K$, such that the representation
$A\rightarrow K$ is irreducible, coincides with $K$, we prove the
following result.

\medskip
\begin{corollary}
\label{c*alg}
The $C^*$--algebra generated by $\Pol(\sph2n)$ is isomorphic to $\tilde K$,
the minimal unitization of compacts.
\end{corollary}
\medskip
\begin{remark}{\rm The
$C^*$--algebra of quantum even spheres is independent of the
classical dimension. This is not as strange as it may appear; the
$C^*$--algebra level usually reflects the topology of the space of
leaves on the underlying Poisson bracket which is the same in all
cases.\fidi}\end{remark}
\begin{remark}{\rm As already hinted in
the introduction, starting from any $C^*$--algebra $A$ with at
least one character $\varepsilon_x$ and from a quantum two-sphere
$B$ with a character $\varepsilon_0$ one could define a
topological double suspension:
\begin{equation}\label{qsosp}
 S^2_qA:=\{f\in A\otimes B\,\big|\, (\varepsilon_x\otimes id)(f)=
 (id\otimes\varepsilon_0)(f)\in{\mathbb C}\}\, .
\end{equation}
It is then a trivial remark that such construction applied to
standard Podle\'s sphere is stable. What is less trivial is the
fact that such algebras \emph{quantize} a whole family of
polynomial quantum even spheres.\fidi}\end{remark}

\medskip
\medskip
Thanks to Corollary (\ref{c*alg}) we can conclude that
$K^0(C(\Sigma^{2n}_q))=\Ze^2$ for each $n$, see \cite{MNW}; each
polynomial sphere $\Pol(\Sigma^{2n}_q)$ will provide different
representatives of the same class in $K$-homology. Let us describe
them explicitly, along the same lines of \cite{MNW}.

The first one $[\epsilon]$ is the pullback by $\epsilon:
C(\Sigma^{2n}_q)\rightarrow \Co$ of the generator of $K^0(\Co)$.
By analogy with the classical case we say that its character
$\epsilon$ computes the rank of the vector bundle.

Let us describe the second and more interesting generator. Let the
Hilbert space be $H=\ell^2(\Na)^{\otimes 2}\oplus
\ell^2(\Na)^{\otimes 2}$ and $\pi=\left (
\begin{array}{rr}\sigma_n & 0\\0
&\epsilon\end{array} \right ),\ F=\left ( \begin{array}{rr}0 &
1\\1 &0\end{array} \right )$. We have that $(H,\pi,F)$ is a
1--summable Fredholm module whose character is the zero cyclic
cocycle $\tr_{\sigma_n}=\tr\; (\sigma_n -\epsilon)$. Let us denote
with $[\tr_{\sigma_n}]$ its class; we say that $\tr_{\sigma_n}$
computes the charge. In particular, from (6) we get
\begin{equation}\label{character}
\tr_{\sigma_n}(1)=0 \,,~~~~
\tr_{\sigma_n}(t)=\frac{1}{(1-q^2)^n}~~~.
\end{equation}

\bigskip
\bigskip
\section{Algebraic projectors and the Chern--Connes \\pairing }

Let us now come to the construction of non trivial quantum vector
bundles on spheres. Since $C(\Sigma^{2n}_q)=\tilde{K}$ we have
that $K_0(C(\Sigma^{2n}_q))=\Ze^2$. In this section we will
introduce two algebraic generators, {\it i.e.} two projectors with
entries in $\Pol(\Sigma^{2n}_q)$ whose classes generate $K_0$. The
first one is the trivial one $[1]$; in order to get the non
trivial generator let us go back to the classical case.

The non trivial generator of $K$-theory for the classical even
sphere can be explicitly written in the following way. Let
$G_{2k}^{2n}\in M_{2^{k}}(\sphere^{2n})$ be defined iteratively by
$$
G^{2n}_{2(k+1)} = \left(\begin{array}{cc}G^{2n}_{2k} &  a^*_{k+1}
\cr a_{k+1} & 1-G^{2n}_{2k}
\end{array}\right) \; ~~~ G^{2n}_0 = 1-t \;.
$$
It is easy to verify that $G_{2n}\equiv G^{2n}_{2n}$ is an
idempotent for $n\geq 0$ defining the vector bundle $E_{2n}$ of
rank $2^{n-1}$ and charge $-1$. This construction has the
following geometrical interpretation. Let
$i:\sphere^{2n}\rightarrow\sphere^{2(n+1)}$ defined by
$i(t,a_1,\ldots,a_n)=(t,a_1,\ldots,a_n,0)$ be an embedding of
$\sphere^{2n}$ in $\sphere^{2(n+1)}$ and let $i^*(E_{2(n+1)})$ the
pullback vector bundle on $\sphere^{2n}$. It is clear that
$i^*(E_{2(n+1)})=E_{2n}\oplus E_{2n}'$, where $E'_{2n}$ is the
conjugated vector bundle on $\sphere^{2n}$ of charge $1$.

In the quantum case not every step of this procedure can be
obviously extended. In fact there is no embedding of
$\Sigma^{2n}_q$ in $\Sigma^{2(n+1)}_q$, i.e. there are no algebra
projections from $\Pol(\Sigma^{2(n+1)}_q)$ to
$\Pol(\Sigma^{2n}_q)$. It is possible to adapt the procedure in
order to produce a rank $2^{n-1}$ idempotent for $\Sigma^{2n}_q$;
but it is convenient to lift the construction to $\Re^{2n+1}_q$, a
deformation of the odd plane where the even sphere lives. Let us
introduce $\Pol(\Re^{2n+1}_q)$ as the algebra generated by
$\{x_i,x_i^*,y=y^*\}_{i=1}^n$ with the relations
$$
x_i\, y=q^2 y\, x_i\;,\quad x_i\, x_j=q^{-1} x_j\, x_i\;, \quad
x_i\, x_j^*=q^{3} x_j^*\, x_i\;,\quad i < j\;,
$$
$$
x_i\, x_i^*=q^2x_i^*\, x_i+q^2(1-q^2)\sum_{\ell<i}x_\ell^*\,
x_\ell + (1-q^2)y^2\;.
$$
It is clear that $\Pol(\Sigma^{2n}_q)=\Pol(\Re^{2n+1}_q)/ I_n$,
where $I_n$ is the ideal generated by $q^2\sum_i^n x_i^*x_i -y +
y^2$ and $a_i=p(x_i), t=p(y)$ if $p$ is the projection map. Let
$\phi_n:\Pol(\Re^{2n+1}_q)\rightarrow \Pol(\Re^{2n+1}_q)$ be the
algebra automorphism defined by $\phi_n(x_i)=q^2x_i$ and
$\phi_n(y)=q^2 y$. The existence of this automorphism, that
doesn't pass to the quotient, is actually the reason for this
lifting.

For each $0\leq k\leq n$ let us define $e^{2n}_{2k}\in
M_{2^k}(\Re^{2n+1}_q)$ using the following recursive formula

\begin{equation}
e^{2n}_{2(k+1)} = \left(\begin{array}{cc}e^{2n}_{2k} & C^{2n}_{2k}
x^*_{k+1} \cr C^{2n}_{2k} x_{k+1} & 1-\phi_n(e^{2n}_{2k})
\end{array}\right) \; ~~~ e^{2n}_0 = 1-y \;,
\end{equation}
where $C^{2n}_{2k}$ is the diagonal complex matrix given by
\begin{equation}
C^{2n}_{2k} = \left(\begin{array}{cc}C^{2n}_{2(k-1)} & 0 \cr0 & q
C^{2n}_{2(k-1)}
\end{array}\right)\; \in M_{2^k}(\Co)  \;, ~~~ C_0^{2n} = q \;.
\end{equation}
\medskip
We need the following Lemma.
\begin{lemma}
For each $\ell$ such that $k<\ell\leq n$ we have that
\begin{equation}\label{lemma_idempotente}
\M^{2n}_{2k,\ell}\equiv e^{2n}_{2k} C^{2n}_{2k} x^*_\ell -
C^{2n}_{2k} x^*_\ell \phi(e^{2n}_{2k}) = 0\;. \end{equation}
\end{lemma}
\noindent {\it Proof}. We prove the result by induction on $k$.
For $k=0$ it is equivalent to $yx_\ell^* = q^2 x_\ell^* y$. Let us
suppose (\ref{lemma_idempotente}) true for $k$ and let us show it
for $k+1$. By direct computation we get
$[\M^{2n}_{2(k+1),\ell}]_{11}= \M^{2n}_{2k,\ell}$,
$[\M^{2n}_{2(k+1),\ell}]_{22}=-q^{-1}\phi_n(\M^{2n}_{2k,\ell})$
and
$[\M^{2n}_{2(k+1),\ell}]_{12}=q(C^{2n}_{2k})^2(x^*_{k+1}x^*_\ell
-q x^*_\ell x^*_{k+1})$. Using the inductive hypothesis and
relations in $\Pol(\Re^{2n+1}_q)$ we conclude that
$\M^{2n}_{2(k+1),\ell}=0$. \fidi

\medskip
We now prove the main result of this section.
\begin{proposition}
For each $k\leq n$ we have
\begin{equation}
\label{intermediate_idempotents} (e^{2n}_{2k})^2 - e ^{2n}_{2k} =
[q^2\sum_{i=1}^k x_i^* x_i - y(1-y) ] q^{-2} (C^{2n}_{2k})^2 \;.
\end{equation}
\end{proposition}
\noindent {\it Proof}. We show the result by induction on $k$. For
$k=0$ it is easy to see that it is true. Let us suppose it true
for $k$. By direct computation, using the inductive hypothesis
and equation (\ref{lemma_idempotente}) we get
\begin{eqnarray*}
[(e^{2n}_{2(k+1)})^2 - e^{2n}_{2(k+1)}]_{11} &=&
[q^2\sum_{i=1}^{k+1} x_i^* x_i - y(1-y) ] q^{-2}
(C^{2n}_{2k})^2\cr [(e^{2n}_{2(k+1)})^2 - e^{2n}_{2(k+1)}]_{22}
&=& (C^{2n}_{2k})^2 x_{k+1}x_{k+1}^* + \phi_n(e^{2n}_{2k})^2 -
\phi_n(e^{2n}_{2k})\cr
&=& (C^{2n}_{2k})^2[x_{k+1}x_{k+1}^* + \sum_{i=1}^kq^4x_i^*x_i -
y(1-q^2y)]\cr
&=& (C^{2n}_{2k})^2[q^2 \sum_{i=1}^{k+1} x_i^*x_i - y(1-y)]\cr
[(e^{2n}_{2(k+1)})^2 - (e^{2n}_{2(k+1)}]_{12} &=& e^{2n}_{2k}
C^{2n}_{2k} x^*_{k+1} - C^{2n}_{2k} x^*_{k+1} \phi(e^{2n}_{2k}) =
0\;.
\end{eqnarray*}
Recalling the iterative definition of $C^{2n}_{2(k+1)}$ we
finally get the result
$$
(e^{2n}_{2(k+1)})^2 - e^{2n}_{2(k+1)} = [q^2\sum_{i=1}^{k+1} x_i^*
x_i - y(1-y) ] q^{-2} (C^{2n}_{2{k+1}})^2\;.
$$\fidi

\medskip

It is then clear that for each $n>0$, $G_{2n}=p(e^{2n}_{2n})\in
M_{2^n}(\Pol(\Sigma^{2n}_q))$ is a projector; let us denote with
$[G_{2n}]$ its class in K-theory (both algebraic and topological).
By using the recursive definition of $e^{2n}_{2k}$ it is easy to
compute the matrix trace of $G_{2n}$. In fact the equation
$$
\Tr (e^{2n}_{2(k+1)}) = 2^k + \Tr(e^{2n}_{2k}-\phi_n(e^{2n}_{2k}))
\;, ~~~ \Tr(e^{2n}_0)=1-y\;,$$ is solved by $\Tr(e^{2n}_{2k})=
2^{k-1} - (1-q^2)^k y $ (for $k\geq 1$), so that we have
\begin{equation}
\label{matrix_trace} \Tr(G_{2n})= 2^{n-1} - (1-q^2)^n t\;.
\end{equation}
By recalling the definition of the Fredholm modules $[\epsilon]$
and $[\tr_{\sigma_n}]$, we compute their Chern--Connes pairing
with $G_{2n}$:
$$
\langle[\epsilon],[G_{2n}]\rangle=2^{n-1}\,,~~~~~~\langle
[\tr_{\sigma_n}],[G_{2n}]\rangle = -1 \;.
$$

\medskip
\begin{corollary}
The projector $G_{2n}$ defines a non trivial class both in
$K_0(\Pol(\Sigma^{2n}_q))$ and $K_0(C(\Sigma^{2n}_q))$.
\end{corollary}

\bigskip
\bigskip
\noindent{\bf Aknowledments}. F.B. wants to thank L.Dabrowski,
G.Landi and E.Hawkins for useful discussions on the subject.

\end{document}